\documentclass[12pt,reqno]{amsart}%
\usepackage{amsfonts,anysize,mathptmx,bm}
\usepackage[noadjust]{cite}
\usepackage[colorlinks,urlcolor=blue,citecolor=blue,linkcolor=blue]{hyperref}
\usepackage{amsmath}
\usepackage{amsfonts}
\usepackage{amssymb}
\usepackage{graphicx}
\marginsize{2cm}{2cm}{1cm}{1cm}
\renewcommand{\vec}[1]{\bm{#1}}
\newtheorem{thm}{Theorem}[section]

\newtheorem{lem}[thm]{Lemma}

\newtheorem{prop}[thm]{Proposition}
\theoremstyle{definition}

\theoremstyle{remark}
\newtheorem{rem}[thm]{Remark}
\newtheorem*{pf}{Proof}
\numberwithin{equation}{section}
\allowdisplaybreaks
\AtEndDocument{\bibliographystyle{num}
\bibliography{C:/Users/lao/Desktop/all}\ }
\begin{document}

\title[quasilinear fourth order elliptic equations]{Solutions for quasilinear fourth order elliptic equations on $\mathbb{R}^{N}$
with sign-changing potential}
\author{Shibo Liu
\and Zhihan Zhao}
\dedicatory{Department of Mathematics, Xiamen University, Xiamen 361005, China}
\thanks{This work was supported by NSFC (11671331) and a special fund from Ministry of Education of China (20180707). It was completed while Shibo Liu was
visiting the Abdus Salam International Centre for Theoretical Physics, he would like to thank ICTP for the hospitality.}

\begin{abstract}
We obtain existence and multiplicity results for quasilinear fourth order
elliptic equations on $\mathbb{R}^{N}$ with sign-changing potential. Our
results generalize some recent results on this problem.

\end{abstract}
\maketitle

\section{Introduction}

In this paper we consider the following fourth order quasilinear elliptic
equations on $\mathbb{R}^{N}$,
\begin{equation}
\left\{
\begin{array}
[c]{ll}%
\Delta^{2}u-\Delta u+V(x)u-\dfrac{1}{2}u\Delta(u^{2})=f(x,u) & \text{in
}\mathbb{R}^{N}\text{,}\\
u\in H^{2}(\mathbb{R}^{N})\text{.} &
\end{array}
\right.  \label{e1}%
\end{equation}
Here $N\leq6$, $V\in C(\mathbb{R}^{N})$ is the potential and $f\in
C(\mathbb{R}^{N}\times\mathbb{R})$ is the nonlinearity, $\Delta^{2}%
=\Delta(\Delta)$ is the biharmonic operator.

Fourth order elliptic problems involving the biharmonic operator have been
studied by many authors. Lazer-McKenna \cite{MR1084570} pointed out that this
type of problems furnish a model to study traveling waves in suspension
bridges. For recent results about semilinear biharmonic equations on bounded
domain $\Omega\subset\mathbb{R}^{N}$, we refer the readers to
\cite{MR2101875,MR2472926,MR2440820,MR2514768}.

To the best of our knowledge, the first result about semilinear biharmonic
equations on $\mathbb{R}^{N}$ seems to be Chabrowski-do \'{O} \cite{MR1894788}%
. In Yin-Wu \cite{MR2735556}, for the following fourth order equation
\begin{equation}
\left\{
\begin{array}
[c]{ll}%
\Delta^{2}u-\Delta u+V(x)u=f(x,u) & \text{in }\mathbb{R}^{N}\text{,}\\
u\in H^{2}(\mathbb{R}^{N}) &
\end{array}
\right.  \label{e2}%
\end{equation}
with Laplacian and linear potential terms, the authors obtained a sequence of
high energy solutions of \eqref{e2} assuming that the potential $V$ has positive
infimum and satisfies condition $\left(  V_{1}\right)  $ given below. The
nonlinearity $f(x,u)$ there is odd and superlinear at infinity. Then, their
results were improved by Ye-Tang \cite{MR2927503}, where the case that
$f(x,u)$ is sublinear was also studied with the aid of a critical point
theorem of Kajikiya \cite{MR2152503}. See also Zhang \emph{et.al.}
\cite{MR3071107} for related results for sublinear problems and
\cite{MR3384039} for the asymptotically linear case.

As far as we know, the first work on fourth order quasilinear elliptic
equations \eqref{e1} on $\mathbb{R}^{N}$ is due to Chen \emph{et.al.}
\cite{MR3276713}, where the positive potential $V$ satisfies condition
$\left(  V_{1}\right)  $ below and the nonlinearity $f(x,u)$ is $4$%
-superlinear in the sense of \eqref{e3}. The results in \cite{MR3276713} were
then extended by Cheng-Tang \cite{MR3591225}. Another recent paper about the
problem \eqref{e1} is Che-Chen \cite{MR3784504}, where the nonlinearity $f(x,u)$
is sublinear so that the variational functional is coercive. By applying the
three critical points theorem from Liu-Su \cite{MR1828101} and the classical Clark
theorem, multiple solutions are obtained.

It is interesting to note that unlike many works in this field, in
\cite{MR3591225} the potential $V$ is allowed to be sign-changing. To deal
with this situation, the authors chose a constant $W_{0}>0$ such that%
\[
W(x)=V(x)+W_{0}\geq1
\]
for all $x\in\mathbb{R}^{N}$. They then considered the equivalent problem with
positive potential $W$:%
\begin{equation}
\left\{
\begin{array}
[c]{ll}%
\Delta^{2}u-\Delta u+W(x)u-\dfrac{1}{2}u\Delta(u^{2})=g(x,u) & \text{in
}\mathbb{R}^{N}\text{,}\\
u\in H^{2}(\mathbb{R}^{N})\text{.} &
\end{array}
\right.  \label{ect}%
\end{equation}
Here among other requirements, the new nonlinearity $g(x,t)=f(x,t)+W_{0}t$
satisfies%
\begin{equation}
\left\vert g(x,t)\right\vert \leq C_{1}\left\vert t\right\vert ^{3}%
+C_{2}\left\vert t\right\vert ^{p-1}\text{,\qquad}\left(  x,t\right)
\in\mathbb{R}^{N}\times\mathbb{R} \label{et}%
\end{equation}
for some $C_{1},C_{2}>0$ and $p\in\left(  4,2_{\ast}\right)  $; where
$2_{\ast}$ is the critical Sobolev exponent for $H^{2}(\mathbb{R}^{N})$,
namely $2_{\ast}=+\infty$ for $N\leq4$ and $2_{\ast}=2N/\left(  N-4\right)  $
for $N>4$. Back to the original problem \eqref{e1}, we have%
\[
\lim_{\left\vert t\right\vert \rightarrow0}\frac{f(x,t)}{t}=\lim_{\left\vert
t\right\vert \rightarrow0}\frac{g(x,t)-W_{0}t}{t}=-W_{0}<0
\]
uniformly in $x\in\mathbb{R}^{N}$. Therefore, the results in Cheng-Tang
\cite{MR3591225} have not resolved the case that $V$ is sign-changing and at
the same time
\begin{equation}
\lim_{\left\vert t\right\vert \rightarrow0}\frac{f(x,t)}{t}=0\text{,\qquad
uniformly in $x\in\mathbb{R}^{N}$.} \label{e4}%
\end{equation}

Motivated by the above observation and our early works
\cite{MR3303004,MR3656292} on Schr\"{o}dinger-Poisson systems, the purpose of
the present paper is to study this issue.

Let $H^{2}=H^{2}(\mathbb{R}^{N})$ be the standard Sobolev space. If $V\in
C(\mathbb{R}^{N})$ is bounded from below, we can choose a constant $m>0$ such
that $\tilde{V}(x)=V(x)+m\geq1$ for $x\in\mathbb{R}^{N}$. On the linear
subspace%
\[
X=\left\{  u\in H^{2}\left\vert \int V(x)u^{2}<\infty\right.  \right\}
\text{,}%
\]
where from now on all integrals are taken over $\mathbb{R}^{N}$ except stated
explicitly, we equip with the inner product%
\[
\left(  u,v\right)  =\int\left(  \Delta u\Delta v+\nabla u\cdot\nabla
v+\tilde{V}(x)uv\right)
\]
and corresponding norm $\left\Vert \cdot\right\Vert $. Then $\left(
X,\Vert\cdot\Vert\right)  $ is a Hilbert space that will be denoted by $X$ for
simplicity. Note that if $V\in C(\mathbb{R}^{N})$ is bounded, then $X$ is
precisely the standard Sobolev space $H^{2}$.

Now we are ready to state our assumptions on $V$ and $f$.

\begin{itemize}
\item[$\left(  V_{1}\right)  $] $V\in C(\mathbb{R}^{N})$ is bounded from
below, $\mu(V^{-1}(-\infty,M])<\infty$ for all $M>0$, where $\mu$ is the
Lebesgue measure on $\mathbb{R}^{N}$.

\item[$\left(  V_{2}\right)  $] $V\in C(\mathbb{R}^{N})$ is a bounded function
such that the quadratic form $\mathfrak{B}:X\rightarrow\mathbb{R}$,%
\begin{equation}
\mathfrak{B}(u)=\frac{1}{2}\int\left(  \left\vert \Delta u\right\vert
^{2}+\left\vert \nabla u\right\vert ^{2}+V(x)u^{2}\right)  \label{eqq}%
\end{equation}
is non-degenerate and the negative space of $\mathfrak{B}$ is finite-dimensional.

\item[$\left(  f_{0}\right)  $] $f\in C(\mathbb{R}^{N}\times\mathbb{R})$ and
there exist $C>0$ and $p\in\left(  4,2_{\ast}\right)  $ such that%
\begin{equation}
\left\vert f(x,t)\right\vert \leq C\left(  \left\vert t\right\vert +\left\vert
t\right\vert ^{p-1}\right)  \text{.}\label{esb}%
\end{equation}

\item[$\left(  f_{1}\right)  $] for $\left(  x,t\right)  \in\mathbb{R}%
^{N}\times\mathbb{R}$ we have $0\leq4F(x,t)\leq tf(x,t)$, moreover, for almost
all $x\in\mathbb{R}^{N}$%
\begin{equation}
\lim_{\left\vert t\right\vert \rightarrow\infty}\frac{F(x,t)}{t^{4}}%
=+\infty\text{,\qquad where }F(x,t)=\int_{0}^{t}f(x,s)\mathrm{d}%
s\text{.}\label{e3}%
\end{equation}

\item[$\left(  f_{2}\right)  $] for any $r>0$, we have%
\[
\lim_{\left\vert x\right\vert \rightarrow\infty}\sup_{0<\left\vert
t\right\vert \leq r}\left\vert \frac{f(x,t)}{t}\right\vert =0\text{.}%
\]

\end{itemize}

If $\left(  V_{1}\right)  $ holds, by a well-known compact embedding
established by Bartsch-Wang \cite{MR1349229}, it can been shown as in Chen
\emph{et.al.}\ \cite[Lemma 2.1]{MR3276713} that for $s\in\lbrack2,2_{\ast})$,
the embedding $X\hookrightarrow L^{s}(\mathbb{R}^{N})$ is compact. The
spectral theory of self-adjoint compact operators implies that the eigenvalue
problem%
\begin{equation}
\Delta^{2}u-\Delta u+V(x)u=\lambda u\text{,\qquad}u\in X\label{ee}%
\end{equation}
possesses a complete sequence of eigenvalues%
\[
-\infty<\lambda_{1}\leq\lambda_{2}\leq\cdots\text{,\qquad}\lambda
_{k}\rightarrow+\infty\text{.}%
\]
Each $\lambda_{k}$ has been repeated in the sequence according to its finite
multiplicity. We denote by $\phi_{k}$ the eigenfunction of $\lambda_{k}$ with
$\left\vert \phi_{k}\right\vert _{2}=1$, where $\left\vert \cdot\right\vert
_{q}$ is the $L^{q}$-norm.

Now we are ready to state the main results of this paper.

\begin{thm}
\label{t1}Suppose $\left(  V_{1}\right)  $, $\left(  f_{0}\right)  $ and
$\left(  f_{1}\right)  $ are satisfied. If \eqref{e4} holds and $0$ is not an
eigenvalue of \eqref{ee}, then \eqref{e1} has a nontrivial solution $u\in X$.
\end{thm}

\begin{thm}
\label{t2}Suppose $\left(  V_{1}\right)  $, $\left(  f_{0}\right)  $ and
$\left(  f_{1}\right)  $ are satisfied. If $f(x,\cdot)$ is odd for all
$x\in\mathbb{R}^{N}$, then \eqref{e1} has a sequence of solutions $\left\{
u_{n}\right\}  $ such that $\Phi(u_{n})\rightarrow+\infty$.
\end{thm}

Since $0$ is not an eigenvalue of \eqref{ee}, we may assume that $0\in\left(
\lambda_{\ell},\lambda_{\ell+1}\right)  $ for some $\ell\geq1$. Of course it
is possible that $0<\lambda_{1}$. In this case from the argument we presented
below it is easy to see that the zero function $u=\vec{0}$ is a local
minimizer of $\Phi$, then as in \cite{MR3591225,MR3276713} the mountain pass
theorem of Ambrosetti-Rabinowitz \cite{MR0370183} can be used to get a nonzero
critical point of $\Phi$, which is a nontrivial solution of \eqref{e1}.
Therefore we will omit this easy situation. Note that if $0\in\left(
\lambda_{\ell},\lambda_{\ell+1}\right)  $ for some $\ell\geq1$, the zero
function is not local minimizer of $\Phi$ anymore, this is the main difference
between our Theorems \ref{t1}, \ref{t2} and the results of
\cite{MR3591225,MR3276713}.

For the case that $V$ satisfies $\left(  V_{2}\right)  $, $X$ is exactly the
standard Sobolev space $H^{2}$, we do not have the compact embedding
$X\hookrightarrow L^{s}(\mathbb{R}^{N})$ for $s\in\lbrack2,2_{\ast})$ any
more. But we still have the following result.

\begin{thm}
\label{t3}Suppose $\left(  V_{2}\right)  $, $\left(  f_{0}\right)  $, $\left(
f_{1}\right)  $ and $\left(  f_{2}\right)  $ are satisfied. If \eqref{e4} holds,
then \eqref{e1} has a nontrivial solution $u\in X$.
\end{thm}

\begin{rem}
\label{rr1}Note that in $(f_{1})$, the limit \eqref{e3} is \emph{point-wise}.
Thus, if $a:\mathbb{R}^{N}\rightarrow\left(  0,\infty\right)  $ is continuous
and decay to zero at infinity, $p\in\left(  4,2_{\ast}\right)  $, then%
\[
f(x,t)=a(x)\left\vert t\right\vert ^{p-2}t
\]
satisfies all our assumptions on $f$ in Theorem \ref{t3}.
\end{rem}

Under the assumptions $\left(  V_{1}\right)  $ or $\left(  V_{2}\right)  $,
and $\left(  f_{0}\right)  $, similar to Chen \emph{et.al.}\ \cite[Lemma
2.2]{MR3276713} we can show that the functional $\Phi:X\rightarrow\mathbb{R}$,%
\[
\Phi(u)=\mathfrak{B}(u)+\frac{1}{2}\int u^{2}\left\vert \nabla u\right\vert
^{2}-\int F(x,u)
\]
is well defined and is of class $C^{1}$. The derivative of $\Phi$ is given by%
\[
\langle\Phi^{\prime}(u),v\rangle=\int\left(  \Delta u\Delta v+\nabla
u\cdot\nabla v+V(x)uv\right)  +\int\left(  uv\left\vert \nabla u\right\vert
^{2}+u^{2}\nabla u\cdot\nabla v\right)  -\int f(x,u)v
\]
for $u,v\in X$. Consequently, critical points of $\Phi$ are weak solutions of
problem \eqref{e1}.

To study the functional $\Phi$, it will be convenient to rewrite the quadratic
part $\mathfrak{B}$ in a simpler form. It is well known that, if $\left(
V_{1}\right)  $ holds and $0$ is not an eigenvalue of \eqref{ee}, or if $\left(
V_{2}\right)  $ holds, then there exists an equialent norm $\left\Vert
\cdot\right\Vert _{V}$ on $X$ such that%
\[
\mathfrak{B}(u)=\frac{1}{2}\left(  \left\Vert u^{+}\right\Vert _{V}%
^{2}-\left\Vert u^{-}\right\Vert _{V}^{2}\right)  \text{,}%
\]
where $u^{\pm}$ is the orthorgonal projection of $u$ on $X^{\pm}$ being
$X^{\pm}$ the positive/negative space of $\mathfrak{B}$. Using this
new norm, $\Phi$ can be rewritten as%
\begin{equation}
\Phi(u)=\frac{1}{2}\left(  \left\Vert u^{+}\right\Vert _{V}^{2}-\left\Vert
u^{-}\right\Vert _{V}^{2}\right)  +\frac{1}{2}\int u^{2}\left\vert \nabla
u\right\vert ^{2}-\int F(x,u)\text{.}\label{enf}%
\end{equation}

The paper is organized as follows. In Section 2 we explain why the usual
linking theorem is not applicable to our problem \eqref{e1}. We will prove
Theorem \ref{t1} by applying Morse theory. Therefore we will introduce some
concepts and results of the theory, verify the assumptions required and then
give the proof of Theorem \ref{t1}. In Section 3 we will prove Theorem
\ref{t2} via the symmetric mountain pass theorem. Finally, after
inverstigating the weak continuity of the functional $u\mapsto\int
u^{2}\left\vert \nabla u\right\vert ^{2}$ on $H^{2}$ and its derivative (see Lemma \ref{l4}), we
use Morse theory again to prove Theorem \ref{t3} in Section 4.

\section{Proof of Theorem \ref{t1}}

In this section and the next section, we always assume that $\left(
V_{1}\right)  $ holds. Consider the quadratic form $\mathfrak{B}$ defined in
\eqref{eqq}. The negative space of $\mathfrak{B}$ is given by%
\[
X^{-}=\operatorname*{span}\left\{  \phi_{1},\ldots,\phi_{\ell}\right\}
\text{.}%
\]
Let $X^{+}$ be the orthogonal complement of $X^{-}$ in $X$, then
$X=X^{-}\oplus X^{+}$. It is well known that there exists a constant $\eta>0$
such that%
\begin{equation}
\pm\mathfrak{B}(u)\geq\eta\Vert u\Vert^{2}\text{,\qquad for }u\in X^{\pm
}\text{.}\label{eq}%
\end{equation}

To find critical points of the functionals with indefinite quadratic part, a
natrual idea is to apply the linking theorem. More precisely, set%
\[
N=\left\{  \left.  u\in X^{+}\right\vert \,\Vert u\Vert=\rho\right\}
\text{,\qquad}M=\left\{  \left.  u\in X^{-}\oplus\mathbb{R}^{+}\phi\right\vert
\,\Vert u\Vert\leq R\right\}  \text{,}%
\]
where $\phi\in X^{+}\backslash0$. If $\Phi$ satisfies the Palais-Smale
$\left(  PS\right)  $ condition and for some $R>\rho>0$,%
\begin{equation}
b=\inf_{N}\Phi>\max_{\partial M}\Phi\text{,}\label{e5}%
\end{equation}
then the linking theorem \cite[Theorem 5.3]{MR845785} gives rise to a nonzero
critical point of $\Phi$. In applications, to verify \eqref{e5} we usually need
to show that $\Phi\leq0$ on $X^{-}$.

However, because the second integral $\int u^{2}\left\vert \nabla u\right\vert
^{2}$ in our functional $\Phi$ is positive for $u\neq0$, it seems impossible
to obtain $\Phi|_{X^{-}}\leq0$ even if we assume $F(x,t)\geq0$ for all
$\left(  x,t\right)  \in\mathbb{R}^{N}\times\mathbb{R}$. Therefore unlike many
other indefinite problems (see e.g. \cite{MR2957647,MR1751952,MR2389415}) the
linking theorem is not applicable for our problem \eqref{e1}.

Fortunately, as in our previous works \cite{MR3303004,MR3656292} on
Schrodinger-Poisson systems, we observe that our functional $\Phi$ has a local
linking at the origin. Therefore, we can apply local linking theory (see \cite{MR1312028,MR802575}) and
infinite dimensional Morse theory (see, e.g., Chang \cite{MR1196690} and
Mawhin-Willem \cite[Chapter 8]{MR982267}) to prove Theorem \ref{t1}. We start
by recalling some concepts and results about Morse theory and local linking.

Let $X$ be a Banach space, $\varphi:X\rightarrow\mathbb{R}$ be a $C^{1}%
$-functional, $u$ is an isolated critical point of $\varphi$ and
$\varphi(u)=c$. Then
\[
C_{i}(\varphi,u):=H_{i}(\varphi_{c},\varphi_{c}\backslash\{0\})\text{,\qquad
}i\in\mathbb{N}=\{0,1,2,\ldots\}\text{,}%
\]
is called the $i$-th critical group of $\varphi$ at $u$, where $\varphi
_{c}:=\varphi^{-1}(-\infty,c]$ and $H_{\ast}$ stands for the singular homology
with coefficients in $\mathbb{Z}$.

If $\varphi$ satisfies the $(PS)$ condition and the critical values of
$\varphi$ are bounded from below by $\alpha$, then following Bartsch-Li
\cite{MR1420790}, we define the $i$-th critical group of $\varphi$ at infinity
by
\[
C_{i}(\varphi,\infty):=H_{i}(X,\varphi_{\alpha})\text{,\qquad}i\in
\mathbb{N}\text{.}%
\]
It is well known that the homology on the right hand side does not depend on
the choice of $\alpha$.

\begin{prop}
[{\cite[Proposition 3.6]{MR1420790}}]\label{ap1}If $\varphi\in C^{1}%
(X,\mathbb{R})$ satisfies $(PS)$ and $C_{\ell}(\varphi,\vec{0})\neq C_{\ell
}(\varphi,\infty)$ for some $\ell\in\mathbb{N}$, then $\varphi$ has a nonzero
critical point.
\end{prop}

\begin{prop}
[{\cite[Theorem 2.1]{MR1110119}}]\label{ap2}Suppose $\varphi\in C^{1}%
(X,\mathbb{R})$ has a \emph{local linking} at
$\vec{0}$ with respect to the decomposition $X=\allowbreak Y\oplus Z$, i.e.,
for some $\varepsilon>0$,
\[%
\begin{array}
[c]{ll}%
\varphi(u)\leq0 & \text{for }u\in Y\cap B_{\varepsilon}\text{,}\\
\varphi(u)>0 & \text{for }u\in(Z\backslash\{\vec{0}\})\cap B_{\varepsilon
}\text{,}%
\end{array}
\]
where $B_{\varepsilon}=\left\{  \left.  u\in X\right\vert \,\Vert u\Vert
\leq\varepsilon\right\}  $. If $\ell=\dim Y<\infty$, then $C_{\ell}%
(\varphi,\vec{0})\neq0$.
\end{prop}

Now, we can begin the investigation of our functional $\Phi$.

\begin{lem}
\label{l1}If $\left(  V_{1}\right)  $ $\left(  f_{0}\right)  $ and \eqref{e4}
hold, $0$ is not an eigenvalue of \eqref{ee}, then $\Phi$ has a local linking at
$\vec{0}$ with respect to the decomposition $X=X^{-}\oplus X^{+}$.
\end{lem}

\begin{pf}
As in the proof of \cite[Lemma 2.2]{MR3276713}, we have for all $u\in X$,%
\begin{align}
\int u^{2}\left\vert \nabla u\right\vert ^{2}  &  \leq\left(  \int\left\vert
u\right\vert ^{6}\right)  ^{1/3}\left(  \int\left\vert \nabla u\right\vert
^{3}\right)  ^{2/3}\nonumber\\
&  \leq\left\vert u\right\vert _{6}^{2}\Vert u\Vert_{W^{1,3}}^{2}\leq S\Vert
u\Vert^{4}\text{,} \label{e6}%
\end{align}
where we have used the continuity of the embeddings $X\hookrightarrow
H^{2}\hookrightarrow W^{1,3}$ and $X\hookrightarrow L^{6}$.

By $\left(  f_{0}\right)  $ and \eqref{e4}, for all $\varepsilon>0$, there
exists $C_{\varepsilon}>0$ such that%
\[
\left\vert F(x,t)\right\vert \leq\varepsilon t^{2}+C_{\varepsilon}\left\vert
t\right\vert ^{p}\text{.}%
\]
Therefore, since $p>4$, using \eqref{e6} we see that as $\Vert u\Vert
\rightarrow0$,%
\[
\int u^{2}\left\vert \nabla u\right\vert ^{2}=o(\Vert u\Vert^{2}%
)\text{,\qquad}\int F(x,u)=o(\Vert u\Vert^{2})\text{.}%
\]
Hence as $\Vert u\Vert\rightarrow0$,%
\[
\Phi(u)=\mathfrak{B}(u)+\frac{1}{2}\int u^{2}\left\vert \nabla u\right\vert
^{2}-\int F(x,u)=\mathfrak{B}(u)+o(\Vert u\Vert^{2})\text{.}%
\]
From this estimate and \eqref{eq}, the conclusion of the lemma follows easily.
\end{pf}

Set $g(x,t)=f(x,t)+mt$, then by $\left(  f_{1}\right)  $,%
\begin{align}
G(x,t)    :=&\int_{0}^{t}g(x,s)\mathrm{d}s=F(x,t)+\frac{m}{2}t^{2}\nonumber\\
  \leq&\,\frac{t}{4}g(x,t)+\frac{\tilde{b}}{4}t^{2}\text{,} \label{e10}%
\end{align}
where $\tilde{b}=b+m>0$. Note that from $\left(  f_{1}\right)  $ we have
$G(x,t)\geq0$ for $\left(  x,t\right)  \in\mathbb{R}^{N}\times\mathbb{R}$ and
\begin{equation}
\lim_{\left\vert t\right\vert \rightarrow\infty}\frac{G(x,t)}{t^{4}}%
=+\infty\text{,\qquad for a.e. }x\in\mathbb{R}^{N}\text{.} \label{exe}%
\end{equation}

The functional $\Phi$ can be rewritten as follows:%
\[
\Phi(u)=\frac{1}{2}\Vert u\Vert^{2}+\frac{1}{2}\int u^{2}\left\vert \nabla
u\right\vert ^{2}-\int G(x,u)\text{,}%
\]
with derivative given by%
\[
\langle\Phi^{\prime}(u),v\rangle=(u,v)+\int\left(  uv\left\vert \nabla
u\right\vert ^{2}+u^{2}\nabla u\cdot\nabla v\right)  -\int g(x,u)v\text{.}%
\]

To apply variational methods, it is important to study the convergence of
asymptotically critical sequences. Based on the compact embedding
$X\hookrightarrow L^{s}$ for $s\in[2,2_*)$, it has been shown in Cheng-Tang \cite[Lemma
2.2]{MR3591225} that if $g$ satisfies \eqref{et} then any bounded $\left(
PS\right)  $ sequence $\left\{  u_{n}\right\}  $ of $\Phi$ has a convergent
subsequence. Under our assumptions, the (new) nonlinearity $g$ only satisfies
the weaker condition%
\[
\left\vert g(x,t)\right\vert \leq C_{1}\left\vert t\right\vert +C_{2}%
\left\vert t\right\vert ^{p-1}\text{,\qquad}\left(  x,t\right)  \in
\mathbb{R}^{N}\times\mathbb{R}\text{.}%
\]
However, this is sufficient to get Eq.\ (2.9) in \cite{MR3591225}. In fact,
because up to a subsequence $u_{n}\rightarrow u$ in $L^{2}$ and in $L^{p}$,
\begin{align*}
\left\vert \int\left(  g(x,u_{n})-g(x,u)\right)  \left(  u_{n}-u\right)
\right\vert  &  \leq\left(  C_{1}+C_{2}\right)  \int\left(  \left\vert
u_{n}\right\vert +\left\vert u\right\vert +\left\vert u_{n}\right\vert
^{p-1}+\left\vert u\right\vert ^{p-1}\right)  \left\vert u_{n}-u\right\vert \\
&  \leq C\left(  \left(  \left\vert u_{n}\right\vert _{2}+\left\vert
u\right\vert _{2}\right)  \left\vert u_{n}-u\right\vert _{2}+\left(
\left\vert u_{n}\right\vert _{p}^{p-1}+\left\vert u\right\vert _{p}%
^{p-1}\right)  \left\vert u_{n}-u\right\vert _{p}\right) \\
&  \rightarrow0\text{,\qquad as }n\rightarrow\infty\text{.}%
\end{align*}
Therefore, under the assumptions of Theorem \ref{t1} (or Theorem \ref{t2}),
any $\left(  PS\right)  $ sequence of our functional $\Phi$ also has a
convergent subsequence.

\begin{lem}
\label{l2}If $\left(  V_{1}\right)  $, $\left(  f_{0}\right)  $ and $\left(
f_{1}\right)  $ hold, then $\Phi$ satisfies the $\left(  PS\right)  $ condition.
\end{lem}

\begin{pf}
Let $\left\{  u_{n}\right\}  \subset X$ be a $\left(  PS\right)  $ sequence of
$\Phi$, that is%
\begin{equation}
c:=\sup_{n}\left\vert \Phi(u_{n})\right\vert <\infty\text{,\qquad}\Phi
^{\prime}(u_{n})\rightarrow0\text{.} \label{e7}%
\end{equation}
By the above remark, it suffices to show that $\left\{  u_{n}\right\}  $ is
bounded. Suppose $\left\{  u_{n}\right\}  $ is unbounded, we may assume $\Vert
u_{n}\Vert\rightarrow\infty$. Then using \eqref{e7} and \eqref{e10} we have%
\begin{align}
4c+\Vert u_{n}\Vert &  \geq4\Phi(u_{n})-\langle\Phi^{\prime}(u_{n}%
),u_{n}\rangle\nonumber\\
&  =\Vert u_{n}\Vert^{2}-\int\left(  4G(x,u_n)-g(x,u_{n})u_{n}\right)
\nonumber\\
&  \geq\Vert u_{n}\Vert^{2}-\tilde{b}\int u_{n}^{2}\text{.} \label{e11}%
\end{align}
Let $v_{n}=\Vert u_{n}\Vert^{-1}u_{n}$. Up to a subsequence, by the compact
embedding $X\hookrightarrow L^{2}(\mathbb{R}^{N})$ we see that%
\[
v_{n}\rightharpoonup v\text{\quad in }X\text{,\qquad}v_{n}\rightarrow
v\text{\quad in }L^{2}(\mathbb{R}^{N})\text{,\qquad}v_{n}(x)\rightarrow
v(x)\text{\quad a.e. }\mathbb{R}^{N}%
\]
for some $v\in X$. Multiplying by $\Vert u_{n}\Vert^{-2}$ to both sides of
\eqref{e11} and letting $n\rightarrow\infty$ yield
\[
\tilde{b}\int v^{2}\geq1\text{.}%
\]
Therefore, $v\neq\vec{0}$.

For $x\in\{v\ne0\}$ we have $\left\vert u_{n}%
(x)\right\vert \rightarrow+\infty$. Hence by \eqref{exe},%
\begin{equation}
\frac{G(x,u_{n}(x))}{\Vert u_{n}\Vert^{4}}=\frac{G(x,u_{n}(x))}{u_{n}^{4}%
(x)}v_{n}^{4}(x)\rightarrow+\infty\text{.} \label{e13}%
\end{equation}
Since $\mu(\{v\ne0\})>0$, by Fatou's lemma we deduce from
\eqref{e13} that%
\begin{equation}
\int\frac{G(x,u_{n})}{\Vert u_{n}\Vert^{4}}\geq\int_{v\neq0}\frac{G(x,u_{n}%
)}{\Vert u_{n}\Vert^{4}}\rightarrow+\infty\text{.} \label{eab}%
\end{equation}
It follows from \eqref{e6} and \eqref{e7} that%
\begin{align*}
o(1)=\frac{\Phi(u_{n})}{\Vert u_{n}\Vert^{4}}  &  =\frac{1}{\Vert u_{n}%
\Vert^{4}}\left(  \frac{1}{2}\Vert u_{n}\Vert^{2}+\frac{1}{2}\int u_{n}%
^{2}\left\vert \nabla u_{n}\right\vert ^{2}-\int G(x,u_{n})\right) \\
&  \leq\frac{1}{2\Vert u_{n}\Vert^{2}}+\frac{S}{2}-\int\frac{G(x,u_{n})}{\Vert
u_{n}\Vert^{4}}\rightarrow-\infty\text{,}%
\end{align*}
this is a contradiction. Therefore $\left\{  u_{n}\right\}  $ is bounded in
$X$.
\end{pf}

To investigate $C_{\ast}(\Phi,\infty)$, using the idea of \cite[Lemma
3.3]{MR3656292} we will prove the following lemma. For this purpose we will
use the equivalent norm $\left\Vert \cdot\right\Vert _{V}$ on $X$ and rewrite
$\Phi$ in the form given in \eqref{enf}. Note that by \eqref{e6}, there exists a
constant $S_{1}>0$ such that%
\begin{equation}
\int u^{2}\left\vert \nabla u\right\vert ^{2}\leq S_{1}\left\Vert u\right\Vert
_{V}^{4}\text{.}\label{e66}%
\end{equation}

\begin{lem}
\label{l3}If $\left(  V_{1}\right)  $, $\left(  f_{0}\right)  $ and $\left(
f_{1}\right)  $ hold, $0$ is not an eigenvalue of \eqref{ee}, then there exists $A>0$ such that if $\Phi(u)\leq-A$,
then%
\[
\left.  \frac{\mathrm{d}}{\mathrm{d}t}\right\vert _{t=1}\Phi(tu)<0\text{.}%
\]

\end{lem}

\begin{pf}
Otherwise, there exists a sequence $\left\{  u_{n}\right\}  \subset X$ such
that $\Phi(u_{n})\leq-n$ but%
\[
\langle\Phi^{\prime}(u_{n}),u_{n}\rangle=\left.  \frac{\mathrm{d}}%
{\mathrm{d}t}\right\vert _{t=1}\Phi(tu_{n})\geq0\text{.}%
\]
Consequently, $\Vert u_{n}\Vert_{V}\rightarrow\infty$ and%
\begin{align}
\Vert u_{n}^{+}\Vert_{V}^{2}-\Vert u_{n}^{-}\Vert_{V}^{2}  &  \leq\left(
\Vert u_{n}^{+}\Vert_{V}^{2}-\Vert u_{n}^{-}\Vert_{V}^{2}\right)  +\int\left[
f(x,u_{n})u_{n}-4F(x,u_{n})\right] \nonumber\\
&  =4\Phi(u_{n})-\langle\Phi^{\prime}(u_{n}),u_{n}\rangle\leq-4n\text{.}
\label{eww}%
\end{align}
Let $v_{n}=\Vert u_{n}\Vert_{V}^{-1}u_{n}$ and $v_{n}^{\pm}$ be the orthogonal
projection of $v_{n}$ on $X^{\pm}$. Then up to a subsequence $v_{n}%
^{-}\rightarrow v^{-}$ for some $v^{-}\in X$, because $\dim X^{-}<\infty$.

If $v^{-}\neq0$, then $v_{n}\rightharpoonup v$ in $X$ for some $v\in
X\backslash\vec{0}$. By assumption $\left(  f_{1}\right)  $ we have%
\[
\frac{f(x,t)t}{t^{4}}\geq\frac{4F(x,t)}{t^{4}}\rightarrow+\infty\text{,\qquad
as }t\rightarrow\infty\text{.}%
\]
Thus, similar to the proof of \eqref{eab}, we obtain%
\begin{equation}
\frac{1}{\Vert u_{n}\Vert_{V}^{4}}\int f(x,u_{n})u_{n}\rightarrow
+\infty\text{.}\label{edc}%
\end{equation}
Now, using \eqref{e66} we have a contradiction%
\begin{align*}
0\leq\frac{\langle\Phi^{\prime}(u_{n}),u_{n}\rangle}{\Vert u_{n}\Vert_{V}^{4}}
&  =\frac{1}{\Vert u_{n}\Vert_{V}^{4}}\left(  \left(  \Vert u_{n}^{+}\Vert
_{V}^{2}-\Vert u_{n}^{-}\Vert_{V}^{2}\right)  +2\int u_{n}^{2}\left\vert
\nabla u_{n}\right\vert ^{2}-\int f(x,u_{n})u_{n}\right)  \\
&  \leq o(1)+2S_{1}-\frac{1}{\Vert u_{n}\Vert_{V}^{4}}\int f(x,u_{n}%
)u_{n}\rightarrow-\infty\text{.}%
\end{align*}

Hence, we must have $v^{-}=\vec{0}$. But $\Vert v_{n}^{+}\Vert_{V}^{2}+\Vert
v_{n}^{-}\Vert_{V}^{2}=1$, we deduce $\Vert v_{n}^{+}\Vert_{V}\rightarrow1$.
Now for large $n$ we have%
\[
\Vert u_{n}^{+}\Vert_{V}=\Vert u_{n}\Vert_{V}\Vert v_{n}^{+}\Vert_{V}\geq\Vert
u_{n}\Vert_{V}\Vert v_{n}^{-}\Vert_{V}=\Vert u_{n}^{-}\Vert_{V}\text{,}%
\]
a contradiction to \eqref{eww}.
\end{pf}

\begin{rem}
\label{rk2}We emphasize that the proof of Lemma \ref{l3} does not depend on
the compactness of the embedding $X\hookrightarrow L^{2}$. Therefore the conclusion of Lemma \ref{l3} remains valid if instead of $(V_1)$, $V$ satisfies $(V_2)$.
\end{rem}

\begin{rem}
\label{rk1}Let $B$ be the unit ball in $X$. Using \eqref{e3} and \eqref{e6}, it is
easy to see that for all $u\in\partial B$,%
\[
\Phi(tu)\rightarrow-\infty\text{,\qquad as }t\rightarrow+\infty\text{.}%
\]
Therefore, as in the proof of \cite[Lemma 3.4]{MR3656292}, for $A>0$ large
enough, using Lemma \ref{l3} we can construct a deformation from $X\backslash
B$ to the level set $\Phi_{-A}=\Phi^{-1}(-\infty,-A]$, and deduce%
\begin{equation}
C_{i}(\Phi,\infty)=H_{i}(X,\Phi_{-A})\cong H_{i}(X,X\backslash
B)=0\text{,\qquad for all $i\in\mathbb{N}$.} \label{e14}%
\end{equation}

\end{rem}

\begin{pf}
[Proof of Theorem \ref{t1}]We have shown that $\Phi$ satisfies $\left(
PS\right)  $ and has a local linking at $\vec{0}$ with respect to the
decomposition $X=X^{-}\oplus X^{+}$. Since $\dim X^{-}=\ell$, Proposition
\ref{ap2} yields $C_{\ell}(\Phi,\vec{0})\neq0$. From \eqref{e14} we see that%
\[
C_{\ell}(\Phi,\vec{0})\neq C_{\ell}(\Phi,\infty)\text{.}%
\]
Therefore by Proposition \ref{ap1} we know that $\Phi$ has a nonzero critical
point $u$, which is a nontrivial solution of the problem \eqref{e1}.
\end{pf}

\section{Proof of Theorem \ref{t2}}

To prove Theorem \ref{t2}, we will apply the following symmetric mountain pass
theorem due to Ambrosetti-Rabinowitz \cite{MR0370183}.

\begin{prop}
[{\cite[Theorem 9.12]{MR765240}}]\label{smt}Let $X$ be an infinite dimensional
Banach space, $\varphi\in C^{1}(X,\mathbb{R})$ be even, satisfies $\left(
PS\right)  $ condition and $\varphi(0)=0$. If $X=Y\oplus Z$ with $\dim
Y<\infty$, and $\varphi$ satisfies

\begin{enumerate}
\item[$\left(  I_{1}\right)  $] there are constants $\rho,\alpha>0$ such that
$\varphi|_{\partial B_{\rho}\cap Z}\geq\alpha$,

\item[$\left(  I_{2}\right)  $] for any finite dimensional subspace $W\subset
X$, there is an $R=R(W)$ such that $\varphi\leq0$ on $W\backslash B_{R(W)}$,
\end{enumerate}
then $\varphi$ has a sequence of critical values $c_{j}\rightarrow+\infty$.
\end{prop}

\begin{lem}
\label{ll}For $i\in\mathbb{N}$, let $Z_{i}=\overline{\operatorname*{span}%
}\left\{  \phi_{i},\phi_{i+1},\ldots\right\}  $ and set%
\[
\beta_{i}=\sup_{u\in Z_{i},\Vert u\Vert=1}\left\vert u\right\vert _{2}\text{.}%
\]
Then $\beta_{i}\rightarrow0$ as $i\rightarrow\infty$.
\end{lem}

\begin{pf}
For $u\in Z_{i}$ with $\Vert u\Vert=1$, we have%
\[
\int\left(  \left\vert \Delta u\right\vert ^{2}+\left\vert \nabla u\right\vert
^{2}+V(x)u^{2}\right)  \geq\lambda_{i}\int u^{2}%
\]
or equivalently (note that $\tilde{V}=V+m$),%
\begin{align*}
1=\Vert u\Vert^{2}  &  =\int\left(  \left\vert \Delta u\right\vert
^{2}+\left\vert \nabla u\right\vert ^{2}+\tilde{V}(x)u^{2}\right) \\
&  \geq\left(  \lambda_{i}+m\right)  \int u^{2}=\left(  \lambda_{i}+m\right)
\left\vert u\right\vert _{2}^{2}\text{.}%
\end{align*}
Therefore, because $\lambda_{i}\rightarrow+\infty$, we have%
\[
\left\vert \beta_{i}\right\vert \leq\frac{1}{\sqrt{\lambda_{i}+m}}%
\rightarrow0\text{.}%
\]

\end{pf}

\begin{pf}
[Proof of Theorem \ref{t2}]Under the assumptions of Theorem \ref{t2}, the
functional $\Phi$ is even and satisfies $\left(  PS\right)  $ condition. It
suffices to show that $\Phi$ verifies the assumptions $\left(  I_{1}\right)  $
and $\left(  I_{2}\right)  $ of Proposition \ref{smt}.

\emph{Verification of }$\left(  I_{1}\right)  $. By $\left(  f_{0}\right)  $,
there exist positive constants $C_{1}$ and $C_{2}$ such that%
\begin{equation}
\left\vert G(x,t)\right\vert \leq C_{1}\left\vert t\right\vert ^{2}%
+C_{2}\left\vert t\right\vert ^{p}\text{.} \label{e15}%
\end{equation}
For $i\in\mathbb{N}$, set $Z_{i}$ and $\beta_{i}$ as in Lemma \ref{ll}. Then
we have $\beta_{i}\rightarrow0$. Choose $k\in\mathbb{N}$ such that%
\[
\lambda=\frac{1}{2}-C_{1}\beta_{k}^{2}>0\text{,}%
\]
then set%
\[
Y=\operatorname*{span}\left\{  \phi_{1},\ldots,\phi_{k-1}\right\}
\text{,\qquad}Z=\overline{\operatorname*{span}}\left\{  \phi_{k},\phi
_{k+1},\ldots\right\}  \text{.}%
\]
We have $X=Y\oplus Z$.

For $u\in Z=Z_{k}$, using \eqref{e15} and note that $p>4$, we have%
\begin{align*}
\Phi(u)  &  =\frac{1}{2}\Vert u\Vert^{2}+\frac{1}{2}\int u^{2}\left\vert
\nabla u\right\vert ^{2}-\int G(x,u)\\
&  \geq\frac{1}{2}\Vert u\Vert^{2}-\int G(x,u) \geq\frac{1}{2}\Vert u\Vert
^{2}-C_{1}\left\vert u\right\vert _{2}^{2}-C_{2}\left\vert u\right\vert
_{p}^{p}\\
&  \geq\left(  \frac{1}{2}-C_{1}\beta_{k}^{2}\right)  \Vert u\Vert^{2}%
-C_{2}C\Vert u\Vert^{p}\\
&  =\lambda\Vert u\Vert^{2}+o(\Vert u\Vert^{2})
\end{align*}
as $\Vert u\Vert\rightarrow0$. From this estimate it is easy to see that
$\left(  I_{1}\right)  $ is valid.

\emph{Verification of }$\left(  I_{2}\right)  $. We only need to show that
$\Phi$ is anti-coercive on any finite dimensional subspace $W$. If this is not
true, there exists $\left\{  u_{n}\right\}  \subset W$ and $A>0$ such that
$\Vert u_{n}\Vert\rightarrow\infty$ but $\Phi(u_{n})\geq-A$. Let $v_{n}=\Vert
u_{n}\Vert^{-1}u_{n}$, then up to a subsequence $v_{n}\rightarrow v$ for some
$v\in W\backslash\vec{0}$, because $\dim W<\infty$. Similar to \eqref{eab} we
deduce that%
\[
\frac{1}{\Vert u_{n}\Vert^{4}}\int G(x,u_{n})\rightarrow+\infty\text{.}%
\]
Thus using \eqref{e6} we have%
\begin{align*}
\Phi(u_{n})  &  =\frac{1}{2}\Vert u_{n}\Vert^{2}+\frac{1}{2}\int u_{n}%
^{2}\left\vert \nabla u_{n}\right\vert ^{2}-\int G(x,u_{n})\\
&  \leq\Vert u_{n}\Vert^{4}\left(  \frac{1}{2\Vert u_{n}\Vert^{2}}+\frac{S}%
{2}-\frac{1}{\Vert u_{n}\Vert^{4}}\int G(x,u_{n})\right)  \rightarrow
-\infty\text{,}%
\end{align*}
a contradiction to $\Phi(u_{n})\geq-A$.
\end{pf}

\begin{rem}
In Cheng-Tang \cite[Theorem 1.3]{MR3591225}, the authors essentially studied
the problem \eqref{ect} with positive potential $W$ and nonlinearity $g(x,u)$
satisfying \eqref{et}. While in our Theorem \ref{t2},
the potential $V$ can be indefinite and the nonlinearity $f(x,u)$ only need to
satisfy the much weaker growth condition \eqref{esb}.
\end{rem}

\section{Proof of Theorem \ref{t3}}

We now assume that $V$ satisfies $\left(  V_{2}\right)  $, then $X$ is exactly
the standard Sobolev space $H^{2}$, the embedding $X\hookrightarrow L^2$ is not compact anymore. Therefore we need to recover the $\left(  PS\right)  $
condition with the help of condition $\left(  f_{3}\right)  $. Firstly, we
need to inverstigate the $C^{1}$-functional $\mathcal{N}:H^{2}\rightarrow
\mathbb{R}$,%
\[
\mathcal{N}(u)=\frac{1}{2}\int u^{2}\left\vert \nabla u\right\vert
^{2}\text{.}%
\]
It is known that the derivative of $\mathcal{N}$ is given by
\[
\langle\mathcal{N}^{\prime}(u),v\rangle=\int\left(  uv\left\vert \nabla
u\right\vert ^{2}+u^{2}\nabla u\cdot\nabla v\right)  \text{,\qquad}u,v\in
H^{2}\text{.}%
\]

\begin{lem}
\label{l4}The functional $\mathcal{N}$ is weakly lower semi-continuous, its
derivative $\mathcal{N}^{\prime}:H^{2}\rightarrow H^{-2}$ is weakly
sequentially continuous.
\end{lem}

\begin{pf}
Let $\left\{  u_{n}\right\}  $ be a sequence in $H^{2}(\mathbb{R}^{N})$ such
that $u_{n}\rightharpoonup u$ in $H^{2}$, we need to show%
\[
\mathcal{N}(u)\leq\varliminf\mathcal{N}(u_{n})\text{,\qquad}\langle
\mathcal{N}^{\prime}(u_{n}),\phi\rangle\rightarrow\langle\mathcal{N}^{\prime
}(u),\phi\rangle\text{.}%
\]
for all $\phi\in H^{2}$.

Since $u_{n}\rightharpoonup u$ in $H^{2}$, by the compactness of the embedding
$H^{2}\hookrightarrow H_{\mathrm{loc}}^{1}$, we have $u_{n}\rightarrow u$ in
$H_{\mathrm{loc}}^{1}$. Therefore up to a subsequence,
\begin{equation}
\nabla u_{n}\rightarrow\nabla u\text{\quad a.e. in }\mathbb{R}^{N}%
\text{,\qquad}u_{n}\rightarrow u\text{\quad a.e. in }\mathbb{R}^{N}\text{.}
\label{e9}%
\end{equation}
By Fatou lemma,%
\[
\mathcal{N}(u)=\int\left\vert \nabla u\right\vert ^{2}u^{2}\leq\varliminf
\int\left\vert \nabla u_{n}\right\vert ^{2}u_{n}^{2}=\varliminf\mathcal{N}%
(u)\text{.}%
\]
Hence $\mathcal{N}$ is weakly lower semi-continuous.

The H\"{o}lder conjugate number of $2_{\ast}$ is%
\[
\left(  2_{\ast}\right)  ^{\prime}=\frac{2_{\ast}}{2_{\ast}-1}=\frac{2N}%
{N+4}\text{.}%
\]
Since $N\leq6$, we have%
\[
r=\frac{N+4}{2N-4}>1\text{,\qquad}r^{\prime}=\frac{r}{r-1}=\frac{N+4}%
{8-N}\text{,\qquad}\frac{2N}{N+4}r^{\prime}=\frac{2N}{8-N}\leq2_{\ast}\text{.}%
\]
For $u\in H^{2}$, since $H^{2}\hookrightarrow W^{1,2^{\ast}}$ and
$H^{2}\hookrightarrow L^{2N/\left(  8-N\right)  }$ continuously, we have%
\begin{align}
\int\left\vert \left\vert \nabla u\right\vert ^{2}u\right\vert ^{\left(
2_{\ast}\right)  ^{\prime}}  &  =\int\left\vert \left\vert \nabla u\right\vert
^{2}u\right\vert ^{2N/\left(  N+4\right)  }=\int\left\vert \nabla u\right\vert
^{4N/\left(  N+4\right)  }\left\vert u\right\vert ^{2N/\left(  N+4\right)
}\nonumber\\
&  \leq\left(  \int\left(  \left\vert \nabla u\right\vert ^{4N/\left(
N+4\right)  }\right)  ^{r}\right)  ^{1/r}\left(  \int\left(  \left\vert
u\right\vert ^{2N/\left(  N+4\right)  }\right)  ^{r^{\prime}}\right)
^{1/r^{\prime}}\nonumber\\
&  =\left(  \int\left\vert \nabla u\right\vert ^{\frac{2N}{N-2}}\right)
^{\left(  2N-4\right)  /\left(  N+4\right)  }\left(  \int\left\vert
u\right\vert ^{2N/\left(  8-N\right)  }\right)  ^{\left(  8-N\right)  /\left(
N+4\right)  }\nonumber\\
&  =\left\vert \nabla u\right\vert _{2^{\ast}}^{4N/\left(  N+4\right)
}\left\vert u\right\vert _{2N/\left(  8-N\right)  }^{2N/\left(  N+4\right)
}\nonumber\\
&  \leq C\Vert u\Vert^{4N/\left(  N+4\right)  }\Vert u\Vert%
^{2N/\left(  N+4\right)  }\text{.} \label{e88}%
\end{align}
Since $\left\{  u_{n}\right\}  $ is bounded in $H^{2}$, from \eqref{e88} we see
that $\left\{  \left\vert \nabla u_{n}\right\vert ^{2}u_{n}\right\}  $ is
bounded in $L^{\left(  2_{\ast}\right)  ^{\prime}}$. From \eqref{e9} we have
$\left\vert \nabla u_{n}\right\vert ^{2}u_{n}\rightarrow\left\vert \nabla
u\right\vert ^{2}u$ a.e.\ in\ $\mathbb{R}^{N}$. The Br\'{e}zis-Lieb lemma
\cite{MR699419} implies that $\left\vert \nabla u_{n}\right\vert ^{2}%
u_{n}\rightharpoonup\left\vert \nabla u\right\vert ^{2}u$ in $L^{\left(
2_{\ast}\right)  ^{\prime}}$. Consequently, for $\phi\in H^{2}$, we have
$\phi\in L^{2_{\ast}}$ and%
\begin{equation}
\int\left\vert \nabla u_{n}\right\vert ^{2}u_{n}\phi\rightarrow\int\left\vert
\nabla u\right\vert ^{2}u\phi\text{.} \label{e16}%
\end{equation}
Similarly, because $N\leq6$ implies $N\leq2_{\ast}$, for $u\in H^{2}$ we have%
\begin{align*}
\int\left\vert u^{2}\nabla u\right\vert ^{\left(  2^{\ast}\right)  ^{\prime}}
&  =\int\left\vert u^{2}\nabla u\right\vert ^{2N/\left(  N+2\right)  }%
=\int\left\vert \nabla u\right\vert ^{2N/\left(  N+2\right)  }\left\vert
u\right\vert ^{4N/\left(  N+2\right)  }\\
&  \leq\left(  \int\left(  \left\vert \nabla u\right\vert ^{2N/\left(
N+2\right)  }\right)  ^{s}\right)  ^{1/s}\left(  \int\left(  \left\vert
u\right\vert ^{4N/\left(  N+2\right)  }\right)  ^{s^{\prime}}\right)
^{1/s^{\prime}}\\
&  =\left(  \int\left\vert \nabla u\right\vert ^{2^{\ast}}\right)  ^{\left(
N-2\right)  /\left(  N+2\right)  }\left(  \int\left\vert u\right\vert
^{N}\right)  ^{4/\left(  N+2\right)  }\\
&  =\left\vert \nabla u\right\vert _{2^{\ast}}^{2N/\left(  N+2\right)
}\left\vert u\right\vert _{N}^{4N/\left(  N+2\right)  }\\
&  \leq C\Vert u\Vert^{2N/\left(  N+2\right)  }\Vert u\Vert
^{4N/\left(  N+2\right)  }\text{,}%
\end{align*}
where%
\[
s=\frac{N+2}{N-2}\text{,\qquad}s^{\prime}=\frac{s}{s-1}=\frac{N+2}{4}\text{.}%
\]
Therefore, for all $i$, the sequence $\left\{  u_{n}^{2}\partial_{i}%
u_{n}\right\}  _{n} $ is bounded in $L^{\left(  2^{\ast}\right)  ^{\prime}}$
and converges point-wise to $u^{2}\partial_{i}u$. Applying the Br\'{e}zis-Lieb
lemma \cite{MR699419} again we deduce%
\[
u_{n}^{2}\nabla u_{n}\rightharpoonup u^{2}\nabla u\text{,\qquad in
}[L^{\left(  2^{\ast}\right)  ^{\prime}}]^{N}\text{.}%
\]
For $\phi\in H^{2}$ (which implies $\phi\in W^{1,2^{\ast}}$ and $\partial
_{i}\phi\in L^{2^{\ast}}$) we have%
\begin{equation}
\int u_{n}^{2}\nabla u_{n}\cdot\nabla\phi\rightarrow\int u^{2}\nabla
u\cdot\nabla\phi\text{.} \label{e17}%
\end{equation}
From \eqref{e16} and \eqref{e17}, for $\phi\in H^{2}$ we have
\begin{align*}
\langle\mathcal{N}^{\prime}(u_{n}),\phi\rangle &  =\int\left\vert \nabla
u_{n}\right\vert ^{2}u_{n}\phi+\int u_{n}^{2}\nabla u_{n}\cdot\nabla\phi\\
&  \rightarrow\int\left\vert \nabla u\right\vert ^{2}u\phi+\int u^{2}\nabla
u\cdot\nabla\phi=\langle\mathcal{N}^{\prime}(u),\phi\rangle\text{.}%
\end{align*}
Therefore, we have proved that $\mathcal{N}^{\prime}$ is weakly sequentially continuous.
\end{pf}

As a consequence of Lemma \ref{l4}, if $u_{n}\rightharpoonup u$ in
$X=H^{2}(\mathbb{R}^{N})$, then
\begin{align}
\varliminf\int\left(  \left\vert \nabla u_{n}\right\vert ^{2}u_{n}\left(
u_{n}-u\right)  +u_{n}^{2}\nabla u_{n}\cdot\nabla\left(  u_{n}-u\right)
\right)   &  =\varliminf\left(  4\mathcal{N}(u_{n})-\langle
\mathcal{N}^{\prime}(u_{n}%
),u\rangle\right)  \nonumber\\
&  \geq4\mathcal{N}(u)-\langle \mathcal{N}^{\prime}(u),u\rangle=0\text{.}\label{e18}%
\end{align}
Having established \eqref{e18}, we can follow the idea of \cite{MR3656292} to
prove that $\Phi$ satisfies the $\left(  PS\right)  $ condition.

\begin{lem}
Under the assumptions of Theorem \ref{t3}, the functional $\Phi$ satisfies the
$\left(  PS\right)  $ condition.
\end{lem}

\begin{pf}
Let $\left\{  u_{n}\right\}  \subset X$ be a $\left(  PS\right)  $ sequence of
$\Phi$. We claim that $\left\{  u_{n}\right\}  $ is bounded. Otherwise, up to
a subsequence we may assume $\Vert u_{n}\Vert_{V}\rightarrow\infty$. Let
$v_{n}=\Vert u_{n}\Vert_{V}^{-1}u_{n}$, then%
\[
v_{n}=v_{n}^{-}+v_{n}^{+}\rightharpoonup v=v^{-}+v^{+}\in X\text{,\qquad}%
v_{n}^{\pm},v^{\pm}\in X^{\pm}\text{.}%
\]
If $v=\vec{0}$, then $v_{n}^{-}\rightarrow v^{-}=\vec{0}$ because $\dim
X^{-}<\infty$. Since%
\[
\Vert v_{n}^{+}\Vert_{V}^{2}+\Vert v_{n}^{-}\Vert_{V}^{2}=1\text{,}%
\]
for $n$ large enough we have%
\[
\Vert v_{n}^{+}\Vert_{V}^{2}-\Vert v_{n}^{-}\Vert_{V}^{2}\geq\frac{1}%
{2}\text{.}%
\]
Therefore, by assumption $\left(  f_{2}\right)  $ we deduce%
\begin{align*}
1+\sup_{n}\left\vert \Phi(u_{n})\right\vert +\Vert u_{n}\Vert_{V} &  \geq
\Phi(u_{n})-\frac{1}{4}\langle\Phi^{\prime}(u_{n}),u_{n}\rangle\\
&  =\frac{1}{4}\Vert u_{n}\Vert_{V}^{2}\left(  \Vert v_{n}^{+}\Vert_{V}%
^{2}-\Vert v_{n}^{-}\Vert_{V}^{2}\right)  +\int\left(  \frac{1}{4}%
f(x,u_{n})u_{n}-F(x,u_{n})\right)  \\
&  \geq\frac{1}{8}\Vert u_{n}\Vert_{V}^{2}\text{,}%
\end{align*}
contradicting $\Vert u_{n}\Vert_{V}\rightarrow\infty$.

If $v\neq\vec{0}$, using $F(x,t)\geq0$ and \eqref{e3}, similar to \eqref{eab}, we have%
\[
\frac{1}{\Vert u_{n}\Vert_{V}^{4}}\int F(x,u_{n})\rightarrow+\infty\text{.}%
\]
This will lead to a contradiction as well because by \eqref{e66}, for $n$ large%
\begin{align*}
\frac{1}{\Vert u_{n}\Vert_{V}^{4}}\int F(x,u_{n}) &  =\frac{\Vert u_{n}%
^{+}\Vert_{V}^{2}-\Vert u_{n}^{-}\Vert_{V}^{2}}{2\Vert u_{n}\Vert_{V}^{4}%
}+\frac{1}{2\Vert u_{n}\Vert_{V}^{4}}\int u_{n}^{2}\left\vert \nabla
u_{n}\right\vert ^{2}-\frac{\Phi(u_{n})}{\Vert u_{n}\Vert_{V}^{4}}\\
&  \leq\frac{S_{1}}{2}+1\text{.}%
\end{align*}
Therefore, $\left\{  u_{n}\right\}  $ is bounded in $X$.

Now, up to a subsequence we have $u_{n}\rightharpoonup u$ in $X$. Hence%
\[
\int\left(  \Delta u_{n}\Delta u+\nabla u_{n}\cdot\nabla u+V(x)u_{n}u\right)
\rightarrow\int\left(  \left\vert \Delta u\right\vert ^{2}+\left\vert \nabla
u\right\vert ^{2}+V(x)u^{2}\right)  =\Vert u^{+}\Vert_{V}^{2}-\Vert
u^{-}\Vert_{V}^{2}\text{.}%
\]
Because $\dim X^{-}<\infty$, we have $u_{n}^{-}\rightarrow u^{-}$ and $\Vert
u_{n}^{-}\Vert_{V}\rightarrow\Vert u^{-}\Vert_{V}$. Consequently%
\begin{align}
o(1) &  =\langle\Phi^{\prime}(u_{n}),u_{n}-u\rangle\nonumber\\
&  =\int\left(  \Delta u_{n}\Delta\left(  u_{n}-u\right)  +\nabla u_{n}%
\cdot\nabla\left(  u_{n}-u\right)  +V(x)u_{n}(u_{n}-u)\right)  \nonumber\\
&  \qquad\qquad+\int\left(  \left\vert \nabla u_{n}\right\vert ^{2}%
u_{n}\left(  u_{n}-u\right)  +u_{n}^{2}\nabla u_{n}\cdot\nabla\left(
u_{n}-u\right)  \right)  -\int f(x,u_{n})\left(  u_{n}-u\right)  \nonumber\\
&  =\left(  \Vert u_{n}^{+}\Vert_{V}^{2}-\Vert u_{n}^{-}\Vert_{V}^{2}\right)
-\left(  \Vert u^{+}\Vert_{V}^{2}-\Vert u^{-}\Vert_{V}^{2}\right)  -\int
f(x,u_{n})\left(  u_{n}-u\right)  \nonumber\\
&  \qquad\qquad+\int\left(  \left\vert \nabla u_{n}\right\vert ^{2}%
u_{n}\left(  u_{n}-u\right)  +u_{n}^{2}\nabla u_{n}\cdot\nabla\left(
u_{n}-u\right)  \right)  \nonumber\\
&  =\Vert u_{n}^{+}\Vert_{V}^{2}-\Vert u^{+}\Vert_{V}^{2}-\int f(x,u_{n}%
)\left(  u_{n}-u\right)  \nonumber\\
&  \qquad\qquad+\int\left(  \left\vert \nabla u_{n}\right\vert ^{2}%
u_{n}\left(  u_{n}-u\right)  +u_{n}^{2}\nabla u_{n}\cdot\nabla\left(
u_{n}-u\right)  \right)  +o(1)\text{.}\label{e19}%
\end{align}
Using the condition $\left(  f_{2}\right)  $, similar to \cite[p.
29]{MR2038142} we can prove%
\[
\varlimsup\int f(x,u_{n})\left(  u_{n}-u\right)  \leq0\text{.}%
\]
We deduce from \eqref{e18} and \eqref{e19} that%
\begin{align*}
&  \varlimsup\left(  \Vert u_{n}^{+}\Vert_{V}^{2}-\Vert u^{+}\Vert_{V}%
^{2}\right)  \\
&  =\varlimsup\left(  \int f(x,u_{n})\left(  u_{n}-u\right)  -\int\left(
\left\vert \nabla u_{n}\right\vert ^{2}u_{n}\left(  u_{n}-u\right)  +u_{n}%
^{2}\nabla u_{n}\cdot\nabla\left(  u_{n}-u\right)  \right)  \right)  \\
&  =\varlimsup\int f(x,u_{n})\left(  u_{n}-u\right)  -\varliminf\int\left(
\left\vert \nabla u_{n}\right\vert ^{2}u_{n}\left(  u_{n}-u\right)  +u_{n}%
^{2}\nabla u_{n}\cdot\nabla\left(  u_{n}-u\right)  \right)  \leq0\text{.}%
\end{align*}
Therefore $\Vert u_{n}^{+}\Vert_{V}\rightarrow\Vert u^{+}\Vert_{V}$ because
\[
\Vert u^{+}\Vert_{V}\le\varliminf\Vert u_{n}^{+}\Vert_{V}\le\varlimsup\Vert u_{n}^{+}\Vert_{V}\le\Vert u^{+}\Vert_{V}\text{.}
\]
We
conclude that $\Vert u_{n}\Vert_{V}\rightarrow\Vert u\Vert_{V}$. Hence from
$u_{n}\rightharpoonup u$ in $X$ we deduce $u_{n}\rightarrow u$ in $X$.
\end{pf}

As we have pointed out in Remark \ref{rk2}, the validity of Lemma \ref{l3}
does not rely on the compactness of the embedding $X\hookrightarrow L^{2}$.
Therefore, using the same proof we deduce that under the assumptions of
Theorem \ref{t3}, there exists $A>0$ such that if $\Phi(u)\leq-A$, then%
\[
\left.  \frac{\mathrm{d}}{\mathrm{d}t}\right\vert _{t=1}\Phi(tu)<0\text{.}%
\]
Applying Remark \ref{rk1}, we deduce $C_{i}(\Phi,\infty)=0$ for all
$i\in\mathbb{N}$. On the other hand, similar to the proof of Lemma \ref{l1} we
can show that $\Phi$ has a local linking at $\vec{0}$ with respect to the
decomposition $X=X^{-}\oplus X^{+}$, thus for $l=\dim X^{-}$ we have
$C_{l}(\Phi,\vec{0})\neq0$. That is
\[
C_{l}(\Phi,\vec{0})\neq C_{l}(\Phi,\infty)\text{.}%
\]
By Proposition \ref{ap1}, $\Phi$ has a nonzero critical point. This completes
the proof of Theorem \ref{t3}.


\begin{thebibliography}{10}

\bibitem{MR0370183}
A.~Ambrosetti, P.~H. Rabinowitz, Dual variational methods in critical point
  theory and applications, J. Functional Analysis 14 (1973) 349--381.

\bibitem{MR1420790}
T.~Bartsch, S.~Li,
  \href{http://dx.doi.org/10.1016/0362-546X(95)00167-T}{Critical point theory
  for asymptotically quadratic functionals and applications to problems with
  resonance}, Nonlinear Anal. 28 (1997) 419--441.

\bibitem{MR2038142}
T.~Bartsch, Z.~Liu, T.~Weth,
  \href{http://dx.doi.org/10.1081/PDE-120028842}{Sign changing solutions of
  superlinear {S}chr\"odinger equations}, Comm. Partial Differential Equations
  29 (2004) 25--42.

\bibitem{MR1349229}
T.~Bartsch, Z.~Q. Wang,
  \href{http://dx.doi.org/10.1080/03605309508821149}{Existence and multiplicity
  results for some superlinear elliptic problems on {${\bf R}^N$}}, Comm.
  Partial Differential Equations 20 (1995) 1725--1741.

\bibitem{MR699419}
H.~Br{\'e}zis, E.~Lieb, \href{http://dx.doi.org/10.2307/2044999}{A relation
  between pointwise convergence of functions and convergence of functionals},
  Proc. Amer. Math. Soc. 88 (1983) 486--490.

\bibitem{MR1894788}
J.~Chabrowski, J.~Marcos~do {\'O},
  \href{http://dx.doi.org/10.1016/S0362-546X(01)00144-4}{On some fourth-order
  semilinear elliptic problems in {$\Bbb R^N$}}, Nonlinear Anal. 49 (2002)
  861--884.

\bibitem{MR1196690}
K.-c. Chang,
  \href{http://dx.doi.org/10.1007/978-1-4612-0385-8}{Infinite-dimensional
  {M}orse theory and multiple solution problems}, Progress in Nonlinear
  Differential Equations and their Applications, 6, Birkh\"auser Boston, Inc.,
  Boston, MA, 1993.

\bibitem{MR3784504}
G.~Che, H.~Chen,
  \href{https://projecteuclid.org/euclid.bbms/1523412051}{Existence of multiple
  nontrivial solutions for a class of quasilinear {S}chr\"odinger equations on
  {$\Bbb R^N$}}, Bull. Belg. Math. Soc. Simon Stevin 25 (2018) 39--53.

\bibitem{MR3303004}
H.~Chen, S.~Liu, \href{http://dx.doi.org/10.1007/s10231-013-0363-5}{Standing
  waves with large frequency for 4-superlinear {S}chr\"odinger-{P}oisson
  systems}, Ann. Mat. Pura Appl. (4) 194 (2015) 43--53.

\bibitem{MR3276713}
S.~Chen, J.~Liu, X.~Wu,
  \href{https://doi.org/10.1016/j.amc.2014.10.021}{Existence and multiplicity
  of nontrivial solutions for a class of modified nonlinear fourth-order
  elliptic equations on {$\Bbb R^N$}}, Appl. Math. Comput. 248 (2014) 593--601.

\bibitem{MR3591225}
B.~Cheng, X.~Tang, \href{http://dx.doi.org/10.1016/j.camwa.2016.10.015}{High
  energy solutions of modified quasilinear fourth-order elliptic equations with
  sign-changing potential}, Comput. Math. Appl. 73 (2017) 27--36.

\bibitem{MR2389415}
Y.~Ding, \href{http://dx.doi.org/10.1142/9789812709639}{Variational methods for
  strongly indefinite problems}, \emph{Interdisciplinary Mathematical
  Sciences}, vol.~7, World Scientific Publishing Co. Pte. Ltd., Hackensack, NJ,
  2007.

\bibitem{MR2152503}
R.~Kajikiya, \href{http://dx.doi.org/10.1016/j.jfa.2005.04.005}{A critical
  point theorem related to the symmetric mountain pass lemma and its
  applications to elliptic equations}, J. Funct. Anal. 225 (2005) 352--370.

\bibitem{MR1751952}
W.~Kryszewski, A.~Szulkin, Generalized linking theorem with an application to a
  semilinear {S}chr\"odinger equation, Adv. Differential Equations 3 (1998)
  441--472.

\bibitem{MR1084570}
A.~C. Lazer, P.~J. McKenna,
  \href{http://dx.doi.org/10.1137/1032120}{Large-amplitude periodic
  oscillations in suspension bridges: some new connections with nonlinear
  analysis}, SIAM Rev. 32 (1990) 537--578.

\bibitem{MR1312028}
S.~J. Li, M.~Willem,
  \href{http://dx.doi.org/10.1006/jmaa.1995.1002}{Applications of local linking
  to critical point theory}, J. Math. Anal. Appl. 189 (1995) 6--32.

\bibitem{MR1828101}
J.~Liu, J.~Su, \href{http://dx.doi.org/10.1006/jmaa.2000.7374}{Remarks on
  multiple nontrivial solutions for quasi-linear resonant problems}, J. Math.
  Anal. Appl. 258 (2001) 209--222.

\bibitem{MR1110119}
J.~Q. Liu, The {M}orse index of a saddle point, Systems Sci. Math. Sci. 2
  (1989) 32--39.

\bibitem{MR802575}
J.~Q. Liu, S.~J. Li, An existence theorem for multiple critical points and its
  application, Kexue Tongbao (Chinese) 29 (1984) 1025--1027.

\bibitem{MR2957647}
S.~Liu, \href{http://dx.doi.org/10.1007/s00526-011-0447-2}{On superlinear
  {S}chr\"odinger equations with periodic potential}, Calc. Var. Partial
  Differential Equations 45 (2012) 1--9.

\bibitem{MR3656292}
S.~Liu, Y.~Wu, \href{https://doi.org/10.1112/blms.12019}{Standing waves for
  4-superlinear {S}chr\"odinger-{P}oisson systems with indefinite potentials},
  Bull. Lond. Math. Soc. 49 (2017) 226--234.

\bibitem{MR982267}
J.~Mawhin, M.~Willem, Critical point theory and {H}amiltonian systems,
  \emph{Applied Mathematical Sciences}, vol.~74, Springer-Verlag, New York,
  1989.

\bibitem{MR765240}
P.~H. Rabinowitz, Minimax methods and their application to partial differential
  equations, in Seminar on nonlinear partial differential equations
  ({B}erkeley, {C}alif., 1983), \emph{Math. Sci. Res. Inst. Publ.}, vol.~2,
  Springer, New York, 1984, 307--320.

\bibitem{MR845785}
P.~H. Rabinowitz, Minimax methods in critical point theory with applications to
  differential equations, \emph{CBMS Regional Conference Series in
  Mathematics}, vol.~65, Published for the Conference Board of the Mathematical
  Sciences, Washington, DC, 1986.

\bibitem{MR2514768}
W.~Wang, A.~Zang, P.~Zhao,
  \href{http://dx.doi.org/10.1016/j.na.2008.10.020}{Multiplicity of solutions
  for a class of fourth elliptic equations}, Nonlinear Anal. 70 (2009)
  4377--4385.

\bibitem{MR2472926}
Y.~Yang, J.~Zhang,
  \href{http://dx.doi.org/10.1016/j.jmaa.2008.08.023}{Existence of solutions
  for some fourth-order nonlinear elliptic problems}, J. Math. Anal. Appl. 351
  (2009) 128--137.

\bibitem{MR2927503}
Y.~Ye, C.-L. Tang,
  \href{http://dx.doi.org/10.1016/j.jmaa.2012.04.041}{Infinitely many solutions
  for fourth-order elliptic equations}, J. Math. Anal. Appl. 394 (2012)
  841--854.

\bibitem{MR2735556}
Y.~Yin, X.~Wu, \href{https://doi.org/10.1016/j.jmaa.2010.10.019}{High energy
  solutions and nontrivial solutions for fourth-order elliptic equations}, J.
  Math. Anal. Appl. 375 (2011) 699--705.

\bibitem{MR2101875}
J.~Zhang, S.~Li, \href{http://dx.doi.org/10.1016/j.na.2004.07.047}{Multiple
  nontrivial solutions for some fourth-order semilinear elliptic problems},
  Nonlinear Anal. 60 (2005) 221--230.

\bibitem{MR3071107}
W.~Zhang, X.~Tang, J.~Zhang,
  \href{https://doi.org/10.1016/j.jmaa.2013.05.044}{Infinitely many solutions
  for fourth-order elliptic equations with general potentials}, J. Math. Anal.
  Appl. 407 (2013) 359--368.

\bibitem{MR3384039}
W.~Zhang, X.~Tang, J.~Zhang,
  \href{https://doi.org/10.1080/00036811.2014.979807}{Ground states for a class
  of asymptotically linear fourth-order elliptic equations}, Appl. Anal. 94
  (2015) 2168--2174.

\bibitem{MR2440820}
J.~Zhou, X.~Wu, \href{https://doi.org/10.1016/j.jmaa.2007.12.020}{Sign-changing
  solutions for some fourth-order nonlinear elliptic problems}, J. Math. Anal.
  Appl. 342 (2008) 542--558.

\end{thebibliography}
\end{document}